\documentclass[leqno,12pt]{article}

\usepackage{amssymb}
\usepackage{euscript}

\evensidemargin\oddsidemargin

\font\smc=cmcsc10 scaled 1200  % Small capitals
\begin{document}
\baselineskip=18pt
\setcounter{page}{1}
\renewcommand{\theequation}{\thesection.\arabic{equation}}  
\newtheorem{theorem}{Theorem}[section] 
\newtheorem{question}[theorem]{Question}
\newtheorem{conjecture}[theorem]{Conjecture}

%%%%%%%%%%%%%%%%%%%%%%%%%%% Environment for remarks
\newenvironment{remark}{
        \vspace{0.3cm} \pagebreak [2] 
        \par 
        \refstepcounter{theorem} 
        \noindent 
        {\bf Remark~\thetheorem\ }}
        {\hfill$\Box$\\[.5mm] } 
%%%%%%%%%%%%%%%%%%%%%%%%%%% 

%%%%%%%%%%%%%%%%%%%%%%%%%%% Equation numberings
\newcommand{\eqnsection}{
\renewcommand{\theequation}{\thesection.\arabic{equation}}
    \makeatletter
    \csname  @addtoreset\endcsname{equation}{section}
    \makeatother}
\eqnsection
%%%%%%%%%%%%%%%%%%%%%%%%%%%

%%%%%%%%%%%%%% Bbb characters
\def\r{{\mathbb R}}        % Real numbers
\def\z{{\mathbb Z}}        % Integers
\def\p{{\mathbb P}}        % Probability
\def\e{{\mathbb E}}        % Expectation
\def\favset{{\mathbb V}}   % Set of favourite sites
\def\fav{V}                % Favourite sites
%%%%%%%%%%%%%%

%%%%%%%%%%%%%%%% Special symbols
\def\ee{\mathrm{e}}                      % Exponential
\def\d{\, \mathrm{d}}                    % Differentiation
\def\ind{{\bf 1}\hskip-2.5pt\mathrm{l}}  % Indicator
%%%%%%%%%%%%%%

%%%%%%%%%%%%%%%%%%%%%%%%%%%%%%%%%%%%%%%%%%%%
%%%%%%%%%%%%%%%%%%%%%%%%%%%%%%%%%%%%%%%%%%%%
%%%%%%%%%%%%%%%%%%%%%%%%%%%%%%%%%%%%%%%%%%%%

%%%%%%%%%%%%%% Beginning of the text

\vglue60pt

\centerline{\LARGE\bf Favourite sites of simple random walk}

\bigskip
\bigskip
\bigskip

\centerline{\smc Zhan Shi and B\'alint T\'oth}

\medskip

\centerline{\it Universit\'e Paris VI \& 
                Technical University Budapest}

\bigskip
\bigskip
\bigskip
\bigskip

{\leftskip=2.5truecm
\rightskip=2.5truecm
\baselineskip=15pt

\noindent{\slshape\bfseries Summary.} 
We survey the current status of the list of questions
related to the favourite (or: most visited) sites of 
simple random walk on $\z$,  
raised by P\'al Erd\H os and P\'al R\'ev\'esz in the 
early eighties. 

\bigskip


\bigskip


} %%%%%% End of narrower

\bigskip
\bigskip
\bigskip
\bigskip

\centerline{\it Dedicated to Professor Endre 
Cs\'aki on the occasion of his 65th birthday}

\bigskip
\bigskip
\bigskip

\section{Introduction}
   \label{S:intro}

Let $(S(n), \, n\in \z_+)$ be a simple symmetric 
random walk on $\z$ with $S(0)=0$. Define

\begin{equation}
    \xi(n,x) := \# \{ 0\le k\le n: \; S(k)=x \} ,
    \label{lt}
\end{equation}

\noindent 
which is referred to as the (site) local time of 
the random walk. For each $n$, consider

\begin{equation}
    \favset (n) := \left\{ x\in \z: \, \xi(n,x) = 
    \max_{y\in \z}
    \xi(n,y) \right\},
    \label{favset}
\end{equation}

\noindent 
which is the set of sites which are the ``most 
visited" by the random walk at step $n$. Following 
Erd\H os and R\'ev\'esz \cite{ER84}, an element of 
$\favset(n)$ is called  ``favourite site". 

The study of $\favset (n)$ was initiated by 
Erd\H os and  R\'ev\'esz \cite{ER84}, and Bass and 
Griffin \cite{Bass-Griffin}, and  received
much research interest from other mathematicians 
since.  Somewhat surprisingly, many of the 
easy-to-formulate and  innocent-looking open
questions raised by Erd\H os and R\'ev\'esz in 
\cite{ER84} and \cite{ER87} remain unanswered so far. 
The present paper aims  not only to have an overview 
upon known results in this field, but also to
insist on unsolved problems in the hope that they will 
meet the interest of, and find solutions from, the reader. 

To illustrate these solved or unsolved problems, we 
mention the following question: what is the probability 
that  $0\in \favset (n)$ for infinitely many $n$? 

Since the random walk is symmetric, one would be tempted 
to think that this probability would be 1. However, the 
correct answer is ``0". In fact, Bass and Griffin 
\cite{Bass-Griffin} proved the following result: 

\begin{equation}
    \lim_{n\to \infty}\, \inf_{x\in\favset(n)} |x| = 
    +\infty, \qquad
    \hbox{\rm a.s.}
    \label{BG}
\end{equation}

\noindent 
In words, the process of favourite sites is transient. 
Actually, Bass and Griffin \cite{Bass-Griffin} showed 
that, for any $\varepsilon>0$, the distance of the set of
favourite sites from the origin goes faster than 
$\sqrt{n}/(\log n)^{11+\varepsilon}$
but slower than $\sqrt{n}/(\log n)^{1-\varepsilon}$,
almost surely. We will discuss this result in more
detail in Section \ref{S:rate}.

Let us mention another innocent-looking question. 
It is trivial that 
$\p\{ \# \favset (n)=1,\, \hbox{i.o.} \}=1$ 
(where ``i.o." stands for infinitely often). A little 
more thinking reveals that 
$\p\{ \# \favset (n)=2, \, \hbox{i.o.} \} =1$. 
Erd\H os and R\'ev\'esz \cite{ER84} asked: 

$$
\p\{ \# \favset (n)\ge 3, \, \hbox{i.o.} \} =?
$$

\noindent 
This problem is still open. However, the following was
recently proved by T\'oth \cite{T00}: 
$\p\{ \# \favset (n)\ge 4, \, \hbox{i.o.} \} =0$. 
For more details, see Subsection \ref{Subs:Q2}.

The rest of the paper is splitted into four sections 
according to the natures of the problems invoked. 
Section \ref{S:rate} concerns the problems of how far 
and how close the favourite sites can be to the
origin. For each of these two problems, we have some 
useful but incomplete information. In Section
\ref{S:ER},  we mention ten questions raised by Erd\H os
and  R\'ev\'esz in \cite{ER84} and \cite{ER87}, and
quoted as Open problems 1--10 in the book of R\'ev\'esz
\cite{R90} (pages 130--131). Only a few of these
questions have found solutions. Section \ref{S:related}
is devoted to some related  problems. More precisely, we
will discuss problems for  rarely visited sites,
favourite edges, and the location  of favourite sites.
Finally in Section \ref{S:others}, we  briefly describe
some known results for favourite sites of other
stochastic processes such as Brownian  motion and
general L\'evy processes. 

\bigskip

\section{Large and small values of the favourite sites}
\label{S:rate}

Throughout this section, we pick  an arbitrary element 
$\fav (n)\in \favset (n)$. 
According to (\ref{BG}), 
$|\fav (n)|$ goes to $\infty$ almost surely 
(when $n\to \infty$). The question is to determine the
rate with which $|\fav (n)|$ goes to infinity. 
The answers are different for lower and upper limits. 

\subsection{Escape rates of favourite sites}
\label{Subs:rate}

As far as the lower limits are concerned, the best 
possible result available so far is due to Bass and 
Griffin \cite{Bass-Griffin}.

\bigskip

\noindent 
{\bf Theorem A ({\rm \cite{Bass-Griffin}}).}
{\it With probability one,}

 \begin{equation}
    \liminf_{n\to \infty} \, 
    {|\fav(n)| \over n^{1/2}(\log n)^{-\gamma}} =
    \cases{0      &if \ \ $\gamma<1$,  \cr 
           \infty &if \ \ $\gamma>11$. \cr}
    \label{BG-rate}
 \end{equation}

\bigskip

Throughout the paper, when we state a limit result 
for $\fav(n)$ as in (\ref{BG-rate}), it is to be 
understood that the convergence holds uniformly in
all $\fav (n)\in \favset (n)$. 

Bass and Griffin \cite{Bass-Griffin} also asked about 
the exact rate of escape of the transient process of 
the favourite sites. This seems to be a very challenging 
problem. Here we formulate it in a weaker form.

\medskip

\begin{question}
    \label{Q:rate}
      Find the value of the constant $\gamma_0$ such 
      that with probability one,

 $$
 \liminf_{n\to \infty} \, 
  {|\fav(n)|\over n^{1/2} (\log n)^{-\gamma}} =
     \cases{0      &{\rm if \ \ }$\gamma<\gamma_0$, \cr 
            \infty &{\rm if \ \ }$\gamma>\gamma_0$. \cr}
 $$

\end{question}

\medskip

According to (\ref{BG-rate}), we must have 
$\gamma_0\in [1, 11]$. 
There is good reason to expect that 
$\gamma_0$ would lie in $[1,2]$. 

\subsection{Upper rates of favourite sites}
\label{Subs:LIL}

If we are interested in the upper limits of 
$|\fav (n)|$, here is an answer in the form of the 
law of the iterated logarithm (LIL), which  was 
discovered independently by Erd\H os and R\'ev\'esz 
\cite{ER84} and Bass and Griffin \cite{Bass-Griffin}.

\bigskip

\noindent 
{\bf Theorem B ({\rm \cite{ER84}, \cite{Bass-Griffin}}).}
{\it With probability one,}

$$
\limsup_{n\to \infty} \; 
{\fav (n)\over (2n \log\log n)^{1/2} }= 1,
\qquad \hbox{\rm a.s.}
$$

\bigskip

Therefore, the favourite site $\fav (n)$ and the 
random walk $S(n)$ satisfy the same LIL. However, a 
closer inspection by Erd\H os and R\'ev\'esz 
\cite{ER84} reveals that they have different
L\'evy upper functions: 

\bigskip

\noindent 
{\bf Theorem C ({\rm \cite{ER84}}).} 
{\it There exists a deterministic sequence 
$(a_n)_{n\ge 1}$ satisfying}

$$
\p \{ S(n) \ge a_n, \, \hbox{\rm i.o.} \} =1, 
\qquad \p\{ \fav(n) \ge
a_n, \, \hbox{\rm i.o.} \} =0 .
$$

\bigskip

The upper functions of $S(n)$ are characterized by the 
Erd\H os--Feller--Kolmogorov--Petrowsky integral test 
(R\'ev\'esz \cite{R90}, p. 35): if $(a_n)_{n\ge 1}$ is a 
non-decreasing sequence of positive numbers, then 

\begin{equation}
    \p\{ S(n) \ge n^{1/2}a_n, \, \hbox{\rm i.o.} \}=
    \cases{ 0 \cr 1 \cr}
      \quad 
    \Longleftrightarrow 
      \quad 
    \sum_n {a_n\over n} \exp\left(
        - {a^2_n\over 2}\right) 
    \cases{ <\infty , \cr =\infty. \cr}
    \label{EFKP}
\end{equation}

\noindent 
Erd\H os and R\'ev\'esz \cite{ER84} and Bass and
Griffin \cite{Bass-Griffin} asked the following
question: 

\medskip

\begin{question}
    \label{Q:test}
      Find an integral test to decide whether 
      $\p\{ \fav (n) \ge a_n,  \, \hbox{\rm i.o.} \} =
      0$. 
\end{question}

\medskip

We believe that a key to Question \ref{Q:test} would 
be to control the upper tail probability of the favourite 
site, formulated here for Brownian motion: let $W$ be a 
standard Brownian motion whose local time process is 
denoted by $(L(t,x); \, t\in \r_+ , \,  x\in \r)$,
i.e., for any bounded Borel function $f$, 

$$
\int_0^t f(W(s))\d s = \int_\r f(x) L(t,x) \d x.
$$

\noindent 
Let $U$ denote the (almost surely) unique favourite 
site at time 1: $L(1,U) = \sup_{x\in \r} L(1,x)$. 
We pose the following

\medskip

\begin{conjecture}
    \label{Q:tail}
      There exists a constant $\nu>1$ such that

 \begin{equation}
      0 < 
      \liminf_{x\to +\infty} \, 
         x^\nu \ee^{x^2/2}\, \p(U>x)
      \le 
      \limsup_{x\to +\infty} \, 
         x^\nu \ee^{x^2/2}\, \p(U>x) 
        < +\infty .
    \label{tail}
 \end{equation}

\end{conjecture}

\medskip

If (\ref{tail}) holds, then we think that we should 
be able to obtain an integral test characterizing the 
upper functions of $\fav (n)$: 

\medskip

\begin{conjecture}
    \label{Q:test-form}
      Let $\nu>1$ be the constant satisfying 
      $(\ref{tail})$. For any non-decreasing sequence 
      $(a_n)_{n\ge 1}$ of positive numbers,

    \begin{equation}
      \p\{ V(n) \ge n^{1/2}a_n, \, \hbox{\rm i.o.} \}=
      \cases{ 0 \cr 1 \cr}
      \quad \Longleftrightarrow \quad 
      \sum_n {a_n^{2-\nu}\over n}
        \exp\left( - {a^2_n\over 2}\right) 
      \cases{ <\infty , \cr =\infty. \cr}    
      \label{test-form}
    \end{equation} 

\end{conjecture}

\medskip

Let us explain why we conjecture $\nu>1$. By the 
trivial inequality $U\le \sup_{0\le t\le 1} W(t)$ and 
the usual estimate for Gaussian tails, it is easily 
seen that if the first inequality in (\ref{tail})
holds, then $\nu\ge 1$. On the other hand, if $\nu$ 
were equal to 1, then (\ref{test-form}) would be the 
same as the Erd\H os--Feller--Kolmogorov--Petrowsky test 
in (\ref{EFKP}) which would contradict Theorem C. This 
leads us to the conjecture $\nu>1$. 

More discussions upon the distribution of $U$ can 
be found in Subsection \ref{Subs:Q5} below.

\bigskip

\section{Ten Erd\H os--R\'ev\'esz questions}
\label{S:ER}

Erd\H os and R\'ev\'esz in \cite{ER84} and \cite{ER87} 
raised many questions about favourite sites of 
random walk, which we quote below. They correspond 
exactly to Questions 1--10 on pp. 130--131 of the book
of R\'ev\'esz \cite{R90}. Recall from (\ref{lt}) that 
$\xi(n,x)$ is the local time of the random walk, and 
from (\ref{favset}) that $\favset (n)$ is the set of 
favourite sites. As before, $\fav (n)$
denotes an arbitrary element of $\favset (n)$. 

\subsection{Joint behaviour of favourite site and
maximum local time}
\label{Subs:Q1}

Let $\xi^*(n) := \sup_{x\in \z} \xi(x,n)$, which is 
the maximum local time. In their proof of Theorem B 
(cf. Subsection \ref{Subs:LIL}), Erd\H os and 
R\' ev\'esz \cite{ER84} noticed that, for any
$\varepsilon>0$, almost surely there exist infinitely 
many $n$ such that simultaneously, 
$\fav (n) \ge (1-\varepsilon)(2n\log\log n)^{1/2}$ and
$\xi^*(n)\ge c(2n\log\log n)^{1/2}$, for some 
positive constant 
$c$ depending on $\varepsilon$. This led them to ask 
the following question: what can be said about the 
joint asymptotics of $\fav (n)$ and $\xi^*(n)$? In
particular, if $\fav (n)$ is close to its maximal 
possible value, how large can $\xi^*(n)$ be? 

If $\fav (n)$ and $\xi^*(n)$ were asymptotically 
independent, then one would expect that the limit 
set of 
$\{ \fav (n)/(2n\log\log n)^{1/2}, 
    \xi^*(n)/(2n\log\log n)^{1/2} \}$ 
should be the half-disc 
$\{ (x,y): \, y\ge0,\, x^2+y^2 \le 1\}$. However, 
the correct answer provided by Cs\'aki et al. 
\cite{CsRS} shows that things do not go exactly 
like this.

\bigskip

\noindent 
{\bf Theorem D ({\rm \cite{CsRS}}).} 
   {\it With probability one, the random sequence  

    $$
      \left(  \, 
       {\fav(n) \over (2n\log\log n)^{1/2}}, \;
      {\xi^*(n) \over (2n\log\log n)^{1/2}}\, 
        \right)_{n\ge 3} 
    $$

    \noindent
    is relatively compact, whose limit set is identical to the 
    simplex $\{ (x,y): \; y\ge 0, \;  |x|+y \le 1 \}$.}

\bigskip

In particular, Theorem D implies Theorem B, and also 
the LIL for the maximum local time of random walk which 
was originally proved by Kesten \cite{Kesten}. 

The proof of Theorem D relies on an invariance
principle for local times due to R\'ev\'esz \cite{R81}
(recalled in (\ref{ltinv}) below) and on the Ray--Knight
theorem for Brownian local time.  

\subsection{Many favourite sites}
\label{Subs:Q2}

At each step (say, $n\to(n+1)$) of the random walk  exactly
one of the following three possibilities occurs:
\begin{itemize}
\item[(1)] the currently occupied site is not
     favourite, $S(n+1)\notin\favset(n+1)$, and thus
     $\favset(n+1)=\favset(n)$ remains unchanged;
\item[(2)] the currently occupied site becomes a new
     favourite besides the favourites of the previous
     time $n$, thus $\favset(n)\subset\favset(n+1)$ and 
     $\favset(n+1)\setminus\favset(n)=\{S(n+1)\}$;
\item[(3)] the random walk revisits a site which was
     already favourite in the previous time $n$, and so
     this new site becomes the only new favourite
     $\favset(n+1)=\{S(n+1)\}\subset\favset(n)$.
\end{itemize}

It follows that the number of favourite sites either remains
unchanged, or increases by one, or drops down to 1. From the
recurrence of the random walk it follows easily that 
$\p\{\#\favset(n)=1,\,\hbox{ i.o.}\}=1$ and 
$\p\{\#\favset(n)=2,\,\hbox{ i.o.}\}=1$. 
The question is 

\begin{equation}
    \p\{ \# \favset (n) = r, \, \hbox{ i.o.}\} = \, ? 
    \qquad r=3,4,\dots
    \label{P=?} 
\end{equation}

\noindent
Erd\H os and R\'ev\'esz \cite{ER84} conjectured that 
this probability should be 0, for any $r\ge3$. 

Recently, T\'oth \cite{T00} proved the following

\bigskip

\noindent 
{\bf Theorem E ({\rm \cite{T00}}).} 
{\it We have,}

$$
\p \{ \# \favset (n) \ge 4, \, \hbox{i.o.}\} =0.
$$

\medskip

The main argument of the proof is based  on some (quite
natural) rearrangements of sums and the Ray--Knight
representation of the local time process stopped at
inverse local times. On the technical level, the proof
of Theorem E relies on controlling the probability
distribution of the number of global maxima of given
height $h\gg1$ of critical Galton--Watson processes.   

\bigskip

\subsection{Frequency of having many favourite sites}
\label{Subs:Q3}

It is intuitively clear that in ``most situations", 
there is only one favourite site. Erd\H os and R\'ev\'esz 
\cite{ER84} were interested in the question about how 
often there are at least two favourite sites. To
formulate their question precisely, let 

$$
\nu_0:=0,
\qquad
\nu_{k+1}:=
\inf\{n>\nu_k:\#\favset(n)>1\},
$$

\noindent
that is: $\nu_k$ is the $k$-th time when there are more than one
favourite sites. Or, alternatively, we can define

$$
\kappa_n:=
\#\{0<j\le n: \#\favset(j)>1\}.
$$

\noindent
Now, the question is: 

$$
\hbox{ Is it true that }
\lim_{k\to\infty}\frac{\nu_k}{k} = \infty,
\hbox{ or, equivalently, }
\lim_{n\to\infty}\frac{\kappa_n}{n} = 0,
\hbox{ with probability one?}
$$

\medskip

The question remains open. 

\subsection{Big jumps of favourite sites}
\label{Subs:Q4}

In order to formulate precisely questions 4, 8 and 9 of the
Erd\H os--R\'ev\'esz list we define the {\it last visited
favourite site} $\ell(n)\in\favset(n)$ inductively, as follows: 

\begin{equation}
\ell(0):=0,
\quad
\ell(n+1) :=
\cases{
\ell(n)  
&if \ \    
$S(n+1)\notin\favset(n+1)$,
\cr
S(n+1)   
&if \ \    
$S(n+1)\in\favset(n+1)$.
\cr
}
\label{lastfav}
\end{equation}

The fourth question concerns how large the jump 
sizes of favourite site can be. The answer is 
formulated in the following LIL. 

\bigskip

\noindent 
{\bf Theorem F ({\rm \cite{CsRS}}).} 
{\it With probability one,}
$$
\limsup_{n\to \infty} \, 
{|\ell(n+1)-\ell(n)| \over (2n\log\log n)^{1/2} }
=1.
$$

\bigskip

In particular, Theorem F tells us that the 
extraordinarily large jumps of favourite site are 
asymptotically comparable to the size of the range 
of the random walk.

The main ingredient in the proof of Theorem F is the
Ray--Knight theorem. 
 
\subsection{Limit law of favourite sites}
\label{Subs:Q5}

The question here is: what is the limit distribution of
$\fav(n)/n^{1/2}$ when $n\to \infty$? 

Here is an argument to show that $\fav(n)/n^{1/2}$ has a 
non-degenerate limit distribution. Indeed, according to 
a theorem of R\'ev\'esz \cite{R81}, possibly in an
enlarged probability space, there exists a coupling for
random walk $(S(n), \, n\in \z_+)$ and Brownian motion 
$(W(t), \, t\in \r_+)$ such that for any
$\varepsilon>0$, 

\begin{equation}
    \sup_{x\in \z} |\xi(n,x) - L(n,x)| = 
    o\left( n^{1/4+ \varepsilon}
    \right) , \qquad \hbox{\rm a.s.},
    \label{ltinv}
\end{equation}

\noindent 
where $\xi$ and $L$ denote the local times of $S(n)$ 
and $W(t)$, respectively. Thus, for any fixed $a>0$, 

\begin{eqnarray*}
    \p\left( {\fav(n) \over n^{1/2}} > a\right) 
    &\le& 
    \p\left( \, \sup_{x> n^{1/2} a}   \xi(n,x) \ge 
                \sup_{x\le n^{1/2} a} \xi(n,x) \right) 
    \cr 
    &\le& 
    \p\left( \, \sup_{x> n^{1/2} a} L(n,x) \ge 
                \sup_{x\le n^{1/2} a} L(n,x) 
                  - n^{1/4+\varepsilon} \right) 
    \cr 
    &=& 
    \p\left( \, \sup_{y> a} L(1,y) \ge 
                \sup_{y\le a} L(1,x) -
                   n^{-1/4+\varepsilon} \right) 
    \cr 
    &\to& 
    \p\left( \, \sup_{y> a} L(1,y) \ge 
                \sup_{y\le a} L(1,y) \right) = 
    \p\left( U>a\right) ,
\end{eqnarray*}

\noindent 
when $n\to \infty$, where $U$ denotes as 
in Subsection \ref{Subs:LIL} the location of the 
maximum (on $\r$) of 
$x\mapsto L(1,x)$: $L(1, U) = \sup_{x\in \r} L(1,x)$.  
Similarly, 

\begin{eqnarray*}
    \p\left( {\fav(n) \over n^{1/2}} > a\right) 
    &\ge&\p\left( \, \sup_{x> n^{1/2} a} \xi(n,x) >
                \sup_{x\le n^{1/2} a} \xi(n,x) \right)
    \cr
    &\ge& 
    \p\left( \, \sup_{x> n^{1/2} a} L(n,x) > 
                 \sup_{x\le n^{1/2} a} L(n,x) 
                   + n^{1/4+\varepsilon} \right) 
    \cr 
    &\to& 
    \p\left( U>a\right) ,
\end{eqnarray*}

\noindent 
so that, 
$$
\lim_{n\to \infty} \, 
{\fav(n) \over n^{1/2}} = U, 
\qquad 
\hbox{\rm in distribution.}
$$

\noindent 
We are grateful to Endre Cs\'aki for having 
communicated to us this simple argument for the weak 
convergence. 

The distribution of $U$ was characterized by 
Theorem 6.2 of Borodin \cite{Borodin}, where a double 
Laplace transform of the limit distribution was computed, 
namely, he got a close-form expression for 
$\e(\ee^{-a \sqrt{A}\, |U|})$ for $a>0$ and an
exponential  variable $A$ which is independent of $U$.
However, the expression for $\e(\ee^{-a \sqrt{A}\,
|U|})$ found in \cite{Borodin} looks very
complicated, involving ratios of Whittaker functions.
We have not  been able to invert the double Laplace
transform, or even to get reasonably good information
for the tail probability of $U$ which would be useful
for the upper functions of $\fav(n)$, see Conjectures
\ref{Q:tail} and \ref{test-form} in Subsection
\ref{Subs:LIL}.  

\subsection{Total number of favourite sites}
\label{Subs:Q6}

Let $\alpha(n)$ denote the number of all the different 
favourite sites up to step $n$, i.e., 
$\alpha(n)= \#(\bigcup_{k=0}^n \favset(k))$. Is
it true that with probability one, for all large $n$, 
$\alpha(n) \le (\log n)^c$ (for some constant $c>0$)? 

The question is still open. Omer Adelman (personal
communication) has a lower bound for $\alpha(n)$, 
proving that with probability one, for all large 
$n$, $\alpha(n) \ge c^* \log n$, for some constant
$c^*>0$.  

\subsection{Durations of favourite sites}
\label{Subs:Q7}

The question is: how long a favourite site can stay 
favourite? More precisely, let 
$\beta(n):= \max\{ j-i: \, 0\le i\le j \le n, \,
\bigcap_{k=i}^j \favset(k) \not= \emptyset\}$. 
In words, $\beta(n)$ is the duration of the longest 
period (before $n$) during which a favourite site 
stays favourite. What can be said about the asymptotic
behaviour of $\beta(n)$? 

No answer available so far. 

\subsection{``Capricious" favourite sites}
\label{Subs:Q8}

If $x$ is a favourite site at some stage, can it 
happen that the favourite site moves away from $x$ but 
later returns to $x$? More precisely, do infinite random  
sequences 
$\dots < c_{n-1}< a_n < b_n < c_n < a_{n+1} < \dots$
of positive integers exist such that 
$\ell(a_n) =\ell(c_n)\not=\ell(b_n)$, for $n=1, 2, \dots$?     
(Recall the definition of the last visited favourite site,
(\ref{lastfav}).) 

The question remains open. 

\subsection{Small jumps of favourite sites}
\label{Subs:Q9}

Consider the random increasing sequence of times, when the last
visited favourite site changes value

$$
\lambda_0:=0,
\quad
\lambda_{k+1}:=
\inf\{n>\lambda_k: \ell(n)\not=\ell(n-1)\} ,
$$

\noindent
and the jump sizes of $\ell(n)$ at these times:

$$
j_k:=|\ell(\lambda_{k+1})-\ell(\lambda_k)|.
$$

\noindent
It seems likely that $j_k \to \infty$ almost surely, with
$k\to\infty$.  How to describe the limit behaviour of $j_k$? 
This question is a companion to the one in Subsection 
\ref{Subs:Q4}. While in Subsection \ref{Subs:Q4} the 
upper behaviour of the jump size was determined with 
satisfying precision, we have not been able to get
any non-trivial information about the lower behaviour. 

\subsection{Occupation times and favourite sites}
\label{Subs:Q10}

The arcsine law says that with big probability, the 
random walk spends a long time on one half of the line 
(say, $\z_+$) and only a short time on the other half 
(in this case $\z_-$). Is it true that the favourite site is
located on the same side where the random walk spends  
the most time? For example, if $(\mu_k)_{k\ge 1}$ is a 
random sequence of integers satisfying 
$\mu_k^{-1} \sum_{j=1}^{\mu_k} \ind_{ \{S_j\ge 0\}} \to
1$ ($k\to \infty$), then Erd\H os and R\'ev\'esz 
\cite{ER84} conjectured that $\fav(\mu_k) \to +\infty$, a.s.  

The answer to the conjecture is no. Relying on the
Ray--Knight theorem and careful analysis of the
sample paths of Bessel processes, Cs\'aki and Shi
\cite{CSprep} prove the existence of a random sequence
$(\mu_k, \, k\ge 1)$ such that $\mu_k^{-1}
\sum_{j=1}^{\mu_k} \ind_{ \{S_j< 0\}} \le  c/(\log\log
\mu_k)^2$ (for some constant $c>0$ and all
$k$), and yet $\fav(\mu_k)<0$ for all $k$. 

\bigskip

\section{Some related questions for simple random walk}
\label{S:related}

\subsection{Rarely visited sites}
\label{Subs:rarely}

\noindent 
{\bf (1)} 
{\it Rarely visited sites}. 
Let $R(n):= \{ S(0), S(1), \cdots, S(n)\}$ be the 
range of the random walk up to step $n$. Erd\H os and 
R\'ev\'esz \cite{ER84} conjectured that for any integer 
$r\ge 3$, the probability that there are infinitely many 
$n$ such that each of the sites of $R(n)$ has been 
visited at least $r$ times up to step $n$ is 0. This was 
disproved by T\'oth \cite{T96}, who showed that for any 
integer $r$, this probability actually equals 1. The
proof of this result again relies on the Ray--Knight
representation of local times. This time one has to
control the probability of the event that the value of
the critical Galton--Watson process drops down to zero
from a given value $r$, uniformly in the initial
condition. 

Another aspect of rarely visited sites was studied by 
Major \cite{Major}. He proved that, if $Z(n)$ denotes 
the number of sites of $R(n)$ which have been visited 
exactly once, then with probability one, 
$\limsup_{n\to \infty} Z(n)/(\log n)^2 = c$, where $c$  
is a finite and positive constant. 

\subsection{Favourite edges}
\label{Subs:edge}

\noindent 
{\bf (2)} 
{\it Favourite edges}. 
Instead of looking at favourite sites (which by 
definition  maximize the {\it site} local time of 
the random walk),  we look at favourite edges which 
maximize the edge local time. We can ask a similar 
question as  Question \ref{Subs:Q2}: are there 
infinitely many  $n$ such that there are at least 3 
favourite edges at  time $n$? T\'oth and Werner \cite{TW} 
proved that with  probability one, there are at most 
finitely many $n$ such that there are at least 4 
favourite edges at time $n$.  This is a simpler predecessor of
Theorem E, cited above. 

\subsection{Location of favourite sites}
\label{Subs:location}

\noindent 
{\bf (3)}
{\it Location of favourite sites}. 
The question is whether it is possible for a favourite 
site to be close to the boundary of the range $R(n)$, 
i.e., close to $\max_{0\le k\le n} S(k)$ or to 
$\min_{0\le k\le n} S(k)$. (The question was originally 
due to Omer Adelman). Cs\'aki and Shi \cite{CsS} proved 
that the answer is no: with probability one, 
$\max_{0\le k\le n} S(k) -\fav(n)$ (and by symmetry, 
$\fav(n)-\min_{0\le k\le n} S(k)$) goes to infinity, 
and the rate of escape was determined. 

\bigskip

\section{Favourite sites of other processes}
\label{S:others}

\noindent 
{\bf (1)} 
{\it Brownian motion}. 
Many  questions and results mentioned above can be 
formulated for the favourite sites of Brownian
motion in the obvious way. Some further discussions 
and questions can be found in Leuridan \cite{Leuridan}. 

\bigskip

\noindent 
{\bf (2)} 
{\it Symmetric stable processes}. 
In Bass et al. \cite{BES}, it was proved that the 
favourite site of a symmetric stable process is also 
transient. Eisenbaum \cite{Eisenbaum} proved a 
collection of interesting results for favourite sites 
of symmetric stable processes. For example, she showed
that for each given $t\ge 0$ there is almost surely a 
unique favourite site at time $t$, and that with 
probability one, for all $t\ge 0$, there are at most 
two favourite sites at time $t$. 

\bigskip

\noindent 
{\bf (3)} 
{\it L\'evy and Markov processes}. 
Most of the results for symmetric stable processes 
can be extended to a larger class of L\'evy and even 
symmetric Markov processes. See Marcus \cite{Marcus}, 
Eisenbaum and Khoshnevisan \cite{Eisenbaum-Khoshnevisan}. 

\bigskip

\noindent 
{\bf (4)} 
{\it Two-dimensional random walk}. 
Let $\fav_2(n)$ denote a favourite site, at time $n$, 
of a simple symmetric random walk on $\z^2$. 
Dembo et al. \cite{DPRZ} proved that 
$(\log \| \fav_2(n)\|)/\log n$ converges to $1/2$ with
probability one, where $\|x\|$ denotes the Euclidean 
modulus in $\r^2$.  

\bigskip

\noindent 
{\bf (5)} 
{\it Transient processes}. 
If $\{ X(t), \; t\ge 0\}$ is a transient process, 
having local time at $t=\infty$, denoted by $L(\infty, x)$, 
then $\{ x\in [0, T]: \, L(\infty,x)=\sup_{y\in [0, T]}
L(\infty, y)\}$ is the set of favourite sites of $X$ 
in $[0,T]$. Bertoin and Marsalle \cite{Bertoin-Marsalle} 
studied the case when $X$ is a Brownian motion with a 
positive drift, and Hu and Shi \cite{HS98} the case when 
$X$ is the modulus of a $d$-dimensional Brownian motion
($d>2$). They obtained respective rates of escape 
(when $T\to \infty$) of favourite sites. 

\bigskip

\noindent 
{\bf (6)} 
{\it Poisson process}. 
Khoshnevisan  and Lewis \cite{Khoshnevisan-Lewis} 
studied favourite  sites of a Poisson process. 
They obtained several laws of  the iterated logarithm. 

\bigskip

\noindent 
{\bf (7)} 
{\it Random walk in random environment}. 
The favourite sites can be defined for nearest neighbour random
walk in random environment, exactly as for the  usual  random
walk.  Hu and Shi  \cite{HS00} considered the recurrent case
and proved that the process of favourite sites is 
again transient, and the escape rate was characterized 
via an integral test. The latter question is still open 
for the usual simple random walk, see Question 
\ref{Q:rate} above. The problem of escape rates of 
favourite sites is the only problem we are aware of, 
which is solved for random walk in random environment, 
but which is open for the usual random walk. 

\bigskip
\bigskip
\bigskip

% \noindent {\Large\bf Acknowledgements}
% \bigskip
% 

\bigskip
\bigskip
\bigskip

{\small

\baselineskip=12pt

\noindent 
\begin{tabular}{lll}
&\hskip5pt  Zhan Shi 
    & \hskip30pt B\'alint T\'oth \\ 
&\hskip5pt  Laboratoire de Probabilit\'es UMR 7599
    & \hskip30pt Institute of Mathematics \\   
&\hskip5pt  Universit\'e Paris VI 
    & \hskip30pt Technical University Budapest \\   
&\hskip5pt  4 place Jussieu
    & \hskip30pt Egry J. u. 1 \\   
&\hskip5pt  F-75252 Paris Cedex 05
    & \hskip30pt H-1111 Budapest \\  
&\hskip5pt  France
    & \hskip30pt Hungary \\ 
&\hskip5pt  {\tt zhan@proba.jussieu.fr} 
    & \hskip30pt {\tt balint@math.bme.hu}
\end{tabular}

}

\end{document}